\documentclass[letterpaper, 10 pt, conference]{ieeeconf}  

\IEEEoverridecommandlockouts                              
\usepackage{graphics} 
\usepackage{epsfig} 
\usepackage{amsmath} 
\usepackage{amssymb} 
\usepackage{cite}
\usepackage{algorithm}
\usepackage{algorithmic}
\usepackage{subfigure}
\usepackage{comment}

\title{\LARGE \bf
	High-performance Uncertainty Quantification in Large-scale Virtual Clinical Trials of Closed-loop Diabetes Treatment
}

\author{Asbj{\o}rn Thode Reenberg, Tobias K. S. Ritschel, Bernd Dammann, John Bagterp J{\o}rgensen
	\thanks{This work was partially funded by the IFD Grand Solution project ADAPT-T2D (9068-00056B). *A. T. Reenberg, T. K. S. Ritschel, Bernd Dammann, and  J. B. J{\o}rgensen are with the Department of Applied Mathematics and Computer Science, Technical University of Denmark, DK-2800 Kgs. Lyngby, Denmark. Corresponding author: J. B. J{\o}rgensen (E-mail: {\tt\small jbjo@dtu.dk}).}
}

\begin{document}
\maketitle
\thispagestyle{empty}
\pagestyle{empty}
\begin{abstract}
	In this paper, we propose a virtual clinical trial for assessing the performance and identifying risks in closed-loop diabetes treatments. Virtual clinical trials enable fast and risk-free tests of many treatment variations for large populations of fictive patients (represented by mathematical models). We use closed-loop Monte Carlo simulation, implemented in high-performance software and hardware, to quantify the uncertainty in treatment performance as well as to compare the performance in different scenarios or of different closed-loop treatments. Our software can be used for testing a wide variety of control strategies ranging from heuristical approaches to nonlinear model predictive control. We present an example of a virtual clinical trial with one million patients over 52~weeks, and we use high-performance software and hardware to conduct the virtual trial in 1~h and 22~min.

\end{abstract}

\section{Introduction}
Clinical trials of medical treatments are crucial to ensuring a high level of safety and efficacy. However, they are also very expensive and time-consuming. Therefore, it is important to assess the treatment performance, identify potential risks, and rigorously compare with state-of-the-art prior to the trials. \emph{Virtual} clinical trials are used for exactly this purpose.
In a virtual clinical trial, each patient is represented by a mathematical model and the clinical trial is simulated using high-performance software and hardware. The simulation is carried out for a large population of virtual patients, many different scenarios, and several variations of the treatment. This allows for thorough and \emph{fast} testing of a large variety of different treatment designs.

In this paper, we specifically consider the treatment of type 1 diabetes~(T1D). One in eleven adults suffer from diabetes (both types), and in 2019, 10\% of the global health expenditure (USD 760 billion) was spent on diabetes~\cite{IDF:2019}. Due to autoimmune destruction of $\beta$-cells, people with T1D are unable to produce insulin. Consequently, life-long treatment involving daily injections of insulin is necessary to avoid elevated blood glucose (BG) levels (hyperglycemia), which can lead to several complications and chronic conditions~\cite{Riddle:etal:2018}.
The BG concentration must be measured in order to determine an appropriate insulin dose. While too little insulin results in hyperglycemia, too much insulin results in hypoglycemia (low BG levels), which can be fatal in very severe cases.

As monitoring the BG and determining the appropriate insulin dose is laborious, there is significant interest in automated closed-loop diabetes treatment systems. Such systems are referred to as artificial pancreases (APs), and they consist of 1)~a continuous glucose monitor (CGM), 2)~a control algorithm that determines the insulin dose, and 3)~a pump which delivers the insulin to the patient (it is possible to use other hormones, e.g., glucagon or amylin, in addition to insulin).
Many AP algorithms have been proposed and the majority is based on heuristics~\cite{Capel:etal:2014}, proportional-integral-derivative~(PID) control~\cite{Marchetti2:etal:2008}, fuzzy logic~\cite{Biester:etal:2019}, linear model predictive control~(MPC)~\cite{Eren-Oruklu:etal:2009, schmidt2013model, Boiroux1:etal:2018}, or nonlinear MPC~\cite{Hovorka:etal:2004, Boiroux2:etal:2018}.

All AP algorithms contain hyperparameters, e.g., the gains in a PID controller, and performance assessment is essential to choosing suitable values for these. There exist both process-independent performance measures for this purpose, e.g., setpoint deviation and variance, and process-specific measures~\cite{Astrom:1989, Harris:1989, Jelali:2006}.
For diabetes treatment, including closed-loop systems, time-in-range (TIR), Hb1Ac values, and the probability of severe hypoglycemia~\cite{Holt:etal:2021, Calsbeek:etal:2013, Lal:etal:2021, Hajizadeh:etal:2019, Sejersen:etal:2021} are commonly used (process-specific) performance measures.
%
The performance of an AP can vary significantly between T1D patients due to differences in physiology (e.g., in pharmacodynamics and pharmacokinetics). Therefore, it is necessary to evaluate the performance measures for a large population of patients in order to accurately estimate the uncertainty. However, due to computational limitations of standard software, it is common to evaluate the performance using only a small number of patients and over a short time span (days or a few weeks).

In this work, we describe an approach for high-performance uncertainty quantification of the performance of AP algorithms in large-scale long-term virtual clinical trials. The approach involves mathematical models based on stochastic differential equations (SDEs), and we use closed-loop Monte Carlo simulation to quantify the performance uncertainty. Furthermore, we propose multiple ways of 1)~visualizing the uncertainty in the performance measures and 2)~comparing the performance for different scenarios or AP algorithms. We implement the Monte Carlo simulation and the AP in parallelized high-performance C code, and the computations are carried out on a high-performance computing (HPC) cluster. Finally, we present a numerical example of a virtual clinical trial with one million patients over 52 weeks which can be carried out in 1~h and 22~min.

The remaining part of this paper is organized as follows. In Section~\ref{sec:VirtualClinicalTrial}, we describe the virtual clinical trial, and in Section~\ref{sec:MonteCarloSimulation}, we present the approach for uncertainty quantification of AP algorithms. Section~\ref{sec:Example:VirtualClinicalTrial} contains the numerical examples, and conclusions are presented in Section~\ref{sec:Conclusion}.

\section{Virtual clinical trial}\label{sec:VirtualClinicalTrial}
The virtual clinical trial consists of 1)~a population of patients, 2)~a protocol containing the trial activities (size and duration of meals, intensity and duration of exercise, etc.), 3)~one or more mathematical models of the patients, 4)~values of the model parameters, and 5)~one or more APs (i.e., control algorithms). Furthermore, it is possible to include both incorrectly announced and unannounced meals and exercise which are some of the key challenges that an AP should be able to address. Finally, the virtual clinical trial allows for \emph{stochastic} mathematical models which can represent unmodeled physiological phenomena, uncertain model parameters, and uncertainty related to meals and exercise, as well as noisy sensor measurements.

\subsection{Patients}
The virtual clinical trial contains one million fictive patients. Each patient is represented by the same information that would be available for a real patient. Specifically, each fictive patient is associated with a unique ID and a set of attributes including first and last name, date and place of birth, sex, height, body weight, and resting heart rate.
The height, the body weight and the resting heart rate are sampled from normal distributions, and the date of birth is sampled from a uniform distribution.
%

\subsection{Protocols}\label{sec:protocols}
A protocol consists of a sequence of model \emph{disturbances}, i.e., uncontrolled inputs to the patient model. Common disturbances are meals and exercise.
Each protocol has an ID and for each disturbance, it contains the disturbance type and size as well as time stamps indicating the beginning and end of the disturbance.

Next, to illustrate the concept, we describe a protocol designed to mimic a Northern European lifestyle in terms of meal times, seasons, work weeks, and the number of vacation weeks and public holidays. Furthermore, the protocol involves a high-carb diet (in particular during winter and autumn), which is challenging for APs.
We divide the year into 4 seasons each consisting of 13 weeks, and we assume 6 weeks of vacation and 10 public holidays, represented as an additional 2 weeks of vacation.
Each season is a different combination of three basis weeks; a \textit{standard week}, an \textit{active week}, and a \textit{vacation week}. Furthermore, each basis week is a different combination of four basis days; a \textit{standard day}, an \textit{active day}, a day with a \textit{movie night}, and a day with a \textit{late night}. Table~\ref{tab:WeeksAndSeasons} shows the compositions of the seasons and the weeks.
\begin{table}[t]
		\caption{Compositions of the seasons and the weeks}
		\label{tab:WeeksAndSeasons}
		\begin{tabular}{p{33.75pt}|ccc}
			\multicolumn{3}{l}{\textbf{Compositions of the seasons}} \\[2pt]
			\hline
			Season & Standard week & Active week & Vacation week \\
			\hline
			Winter & 6 & 4 & 3 \\
			Spring & 6 & 6 & 1 \\
			Summer & 7 & 3 & 3 \\
			Autumn & 9 & 3 & 1 \\ \hline
		\end{tabular}

		\vspace{7pt}

		\begin{tabular}{p{33.75pt}|cccc}
			\multicolumn{3}{l}{\textbf{Compositions of the weeks}} \\[2pt]
			\hline
			Week type & Standard day & Active day & Movie night & Late night \\
			\hline
			Standard & 4 & 1 & 1 & 1  \\
			Active   & 3 & 3 & 1 & 0  \\
			Vacation & 5 & 0 & 0 & 2  \\ \hline
		\end{tabular}
\end{table}

\begin{figure}[tb]
	\centering
	\includegraphics[width=0.45\textwidth]{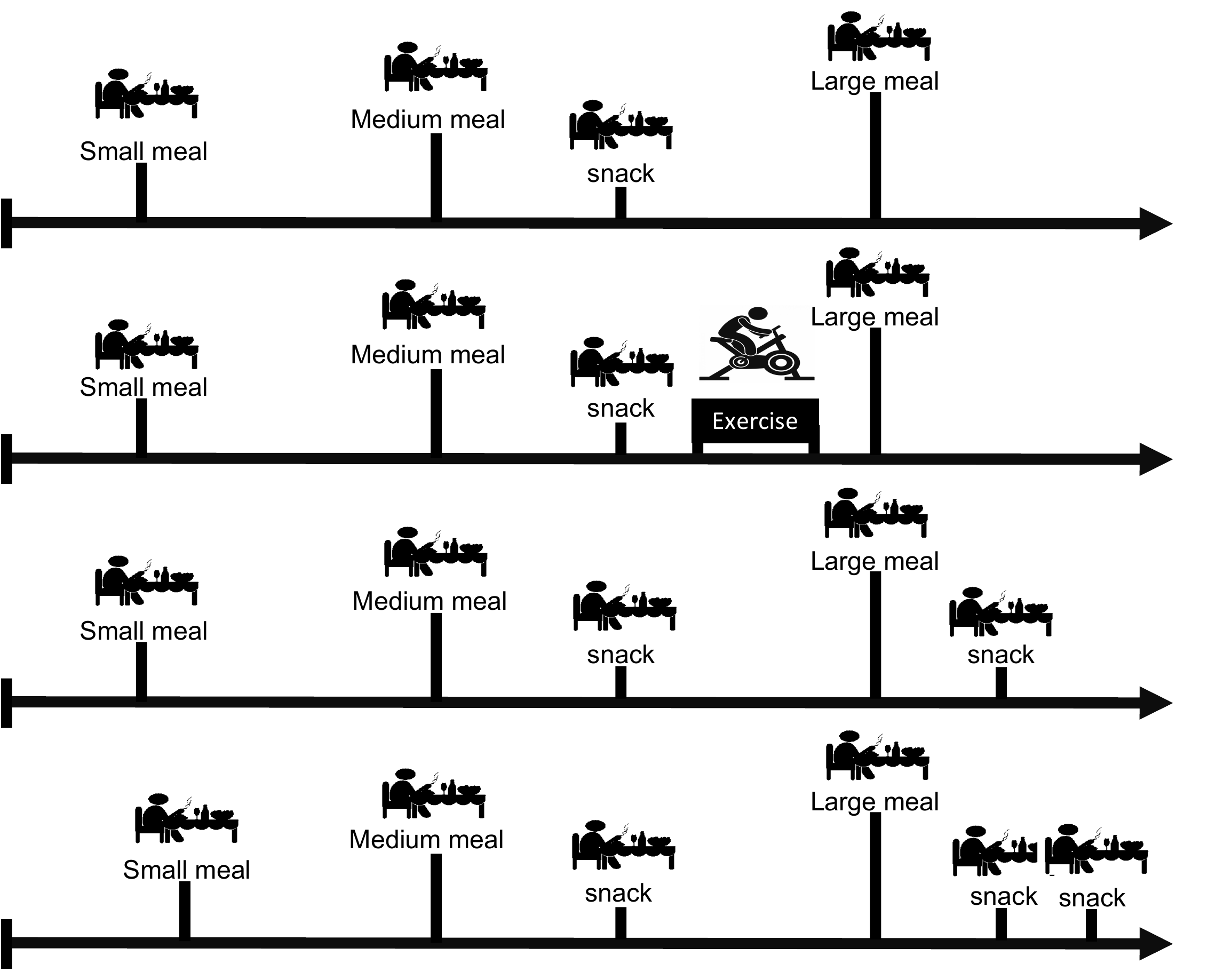}
	\caption{Overview of the basis days in winter and autumn. From the top, we show 1) the \textit{standard day}, 2) the \textit{active day}, 3) the \textit{movie night}, and 4) the \textit{late night}. During the summer and spring, the dinner is a \textit{medium meal} and the snack is before lunch.}
	\label{fig:WinterSummerDays}
\end{figure}

%
%
The patients are less active and eat more during vacation weeks, they also eat more during winter and autumn, and active weeks contain more active days. Compared to the standard day, 1)~the active day has an exercise session, 2)~the movie night has an additional snack in the evening, and 3)~the late night has two additional snacks in the evening.
Fig.~\ref{fig:WinterSummerDays} shows schematics of the basis days, and Table~\ref{tab:MealSizes} shows the meal sizes which depend on the body weight.

\begin{table}[b]
	\begin{center}
		\caption{Weight-dependent meal sizes}
		\label{tab:MealSizes}
		\begin{tabular}{l|cc}
			\hline
			Meal size 	& Amount of carbohydrates & For a 70~kg patient \\
			\hline
			Large meal  & 1.29~g CHO/kg & 90 g CHO \\
			Medium meal & 0.86~g CHO/kg & 60 g CHO \\
			Small meal  & 0.57~g CHO/kg & 40 g CHO \\
			Snack       & 0.29~g CHO/kg & 20 g CHO \\
			\hline
		\end{tabular}
	\end{center}
\end{table}

\subsection{Database}
A key component of the virtual clinical trial is a database containing the fictive patients, protocols, model parameters, and simulation results. We use the open-source database system PostgreSQL. Using a database makes it straightforward to share the virtual clinical trial and compare performance results for the exact same patients and protocols. It also allows the user to carry out the clinical trial for a specific demographic, e.g., people with a certain body weight or sex.
Furthermore, the database includes several basis days and weeks which the user can combine to form new protocols.
Finally, the database can be extended with a graphical user interface in order to 1)~visualize the AP performance and characteristics of the patients and protocols and 2)~add new elements such as patients, protocols, and mathematical models to the database.

\section{Closed-loop Monte Carlo simulation}\label{sec:MonteCarloSimulation}
In Monte Carlo simulation, the uncertainty of a quantity of interest (e.g., the TIR for a T1D patient) is estimated by simulating the system (i.e., the clinical trial) with different values of the uncertain quantities (e.g., the model parameters). In this work, we extend the high-performance toolbox for closed-loop Monte Carlo simulation developed by Wahlgreen et al.~\cite{Wahlgreen:etal:2021} with a low-memory implementation which circumvents the high storage requirements associated with large numbers of long-term simulations.

\subsection{Mathematical model}
The virtual clinical trial can be used for mathematical models of patient physiology in the general stochastic form
\begin{subequations}\label{eq:sde}
	\begin{align}
		\label{eq:sde:initial}
		x(t_0) &= x_0, \\
		\label{eq:sde:state}
		d x(t) &= f(t, x(t), u(t), d(t), p) dt \nonumber \\
		& \quad + \sigma(t, x(t), u(t), d(t), p) dw(t), \\
		\label{eq:sde:output}
		z(t) &= h(t, x(t), p) \\
		\label{eq:sde:observed}
		y(t_k) &= g(t_k, x(t_k), p) + v(t_k).
	\end{align}
\end{subequations}
Here, 
$t$ is time, and the virtual clinical trial starts at time~$t_0$. The states, $x$, represent the physiological state of the patient, e.g., the BG concentration and amount of insulin in the body, and $x_0$ are the initial states. The manipulated inputs, $u$, are the quantities computed by the AP, e.g., the insulin flow rate. The disturbance variables, $d$, represent the uncontrolled inputs, e.g., meals and exercise, and $p$ are model parameters.
The first term in~\eqref{eq:sde:state} is the deterministic drift term and the second term is the stochastic diffusion term.

The AP receives measurements of the observed variables, $y$, obtained from the measurement function, $g$, at discrete points in time, $t_k$, and the outputs, $z$, obtained from the output function, $h$, are the quantities relevant to the control objective of the AP. The standard Wiener process $w$ is used to represent uncertainty. Its increment is distributed as $dw(t) \sim N_{iid}(0,I dt)$, and the measurement noise is assumed to be normally distributed: $v(t_k) \sim N_{iid}(0, R(t_k))$. Furthermore, the inputs are assumed to be piecewise constant between sampling times:
\begin{subequations}
	\begin{align}
		u(t) &= u_k, & t &\in[t_k, t_{k+1}[, \\
		d(t) &= d_k, & t &\in[t_k, t_{k+1}[.
	\end{align}
\end{subequations}
Finally, we stress that the form~\eqref{eq:sde:state} also includes \emph{deterministic} dynamical systems, i.e., ordinary differential equations (ODEs), where $\sigma$ is zero.

\subsection{Control algorithm}
At time $t_k$, 1)~the control state, $x_k^c$, is updated, and 2)~the AP (i.e., the closed-loop feedback control strategy) computes values of the manipulated inputs based on the previous control state and the measurements, $y_k = y(t_k)$:
\begin{subequations}
	\begin{align}
		x_{k+1}^c 	&= \kappa_k(x_k^c, y_k, \bar u_k, \bar y_k, \hat d_k, p_\kappa), \\
		u_k			&= \lambda_k(x_k^c, y_k, \bar u_k, \bar y_k, \hat d_k, p_\mu).
	\end{align}
\end{subequations}
Here, $\bar u_k$ and $\bar y_k$ are setpoints, and $\hat d_k$ are estimates of the disturbances.
This form can represent many types of closed-loop control strategies including heuristic strategies based on physiological insight, PID-based strategies, and MPC (including state estimation). Many different values of the hyperparameters $p_\mu$ and $p_\kappa$ can be tested using the virtual clinical trial.

\subsection{Software and hardware}
The closed-loop Monte Carlo simulation, the mathematical models, and the AP are implemented using high-performance C code which we parallelize for shared-memory architectures using OpenMP. Whenever a simulation is completed, we immediately compute its contribution to performance indicators such as TIR as well as mean, minimum, and maximum BG concentration as functions of time. Subsequently, the simulation is only stored if it is worse than previous simulations according to some criterion (e.g., lowest BG concentration reached).

We use two AMD EPYC 7542 32-core processors with a clock speed of 2.9~GHz \cite{DTU_DCC_resource}. As the Monte Carlo simulation is highly parallelizable, the speedup in computational performance increases almost linearly with the number of cores. Consequently, the parallel implementation runs almost 64 times faster than a corresponding sequential implementation. Furthermore, it would be computationally infeasible to carry out large-scale long-term virtual clinical trials using a sequential Matlab implementation~\cite{Wahlgreen:etal:2021} or similar.

\section{Example of a virtual clinical trial}\label{sec:Example:VirtualClinicalTrial}
In this section, we present an example of a virtual clinical trial involving a million virtual patients following the example protocol described in Section~\ref{sec:protocols}. We briefly describe the mathematical model of the patients' physiology and the AP used in the trial, and we demonstrate how to visualize the uncertainty in the AP performance. Furthermore, we also show how to compare the performance of the AP in two different scenarios; one where the basal rate is correct (trial A) and another where it is underestimated by 50\% (trial B). Finally, the computation time is 1~h and 22~min (with a time step size of 30~s in the simulations). Consequently, several virtual clinical trials can be carried out in a single day.

\subsection{Performance measures}
We divide the BG concentration into 5 ranges~\cite{Holt:etal:2021} given in mmol/L.
Red: severe hypoglycemia (below $3$). Light red: hypoglycemia ($3$--$3.9$). Green: normoglycemia ($3.9$--$10$). Yellow: hyperglycemia ($10$--$13.9$). Orange: severe hyperglycemia (above $13.9$). Furthermore, we also consider the distributions of the total daily basal and bolus insulin as well as bolus glucagon.

\subsection{Patient model}
We use an extension of the model presented by Hovorka et al.~\cite{Hovorka:etal:2004} to represent the pharmacokinetic and pharmacodynamical responses of the virtual patients to carbohydrate absorption and subcutaneous infusion of insulin and glucagon. The model is extended with 1)~a one-state model of the measured BG concentration, 2)~a two-state pharmacokinetic model of subcutaneous glucagon injection~\cite{Wendt:etal:2016}, and 3)~a three-state model of the effect of exercise on the plasma BG concentration~\cite{Rashid:etal:2019}.

The model parameters related to the measured BG concentration, glucagon infusion, and exercise have the same values for all virtual patients. We sample the remaining model parameters from the distributions presented by Hovorka et al.~\cite{Hovorka:etal:2002}. We only use parameter sets where the parameter values are nonnegative, the normally distributed parameters are within one standard deviation from the mean, and the insulin basal rate is at least 0.4~U/h.

\subsection{Artificial pancreas}
We demonstrate the capabilities of the virtual clinical trial using a dual-hormone AP which switches between an insulin mode and a glucagon mode. Glucagon is used to mitigate hypoglycemia. The insulin mode involves 1)~microadjustments of the basal rate, 2)~a meal bolus calculator, 3)~superboli, and 4)~an insulin-to-carb ratio estimator. In the glucagon mode, only microboli are administered. In both modes, a 100~$\mu$g glucagon bolus is administered at the beginning of exercise if the blood glucose concentration is below 7~mmol/L. No insulin is administered in the glucagon mode, and no glucagon, apart from the exercise bolus, is administered in the insulin mode.
The AP uses filtered estimates of the glucose concentration obtained with a low-pass filter, and several hyperparameters have different values depending on whether the patient is exercising or not.

\begin{figure}[tb]
	\centering
	\includegraphics[trim=0 100 0 100, clip,width=0.5\textwidth]{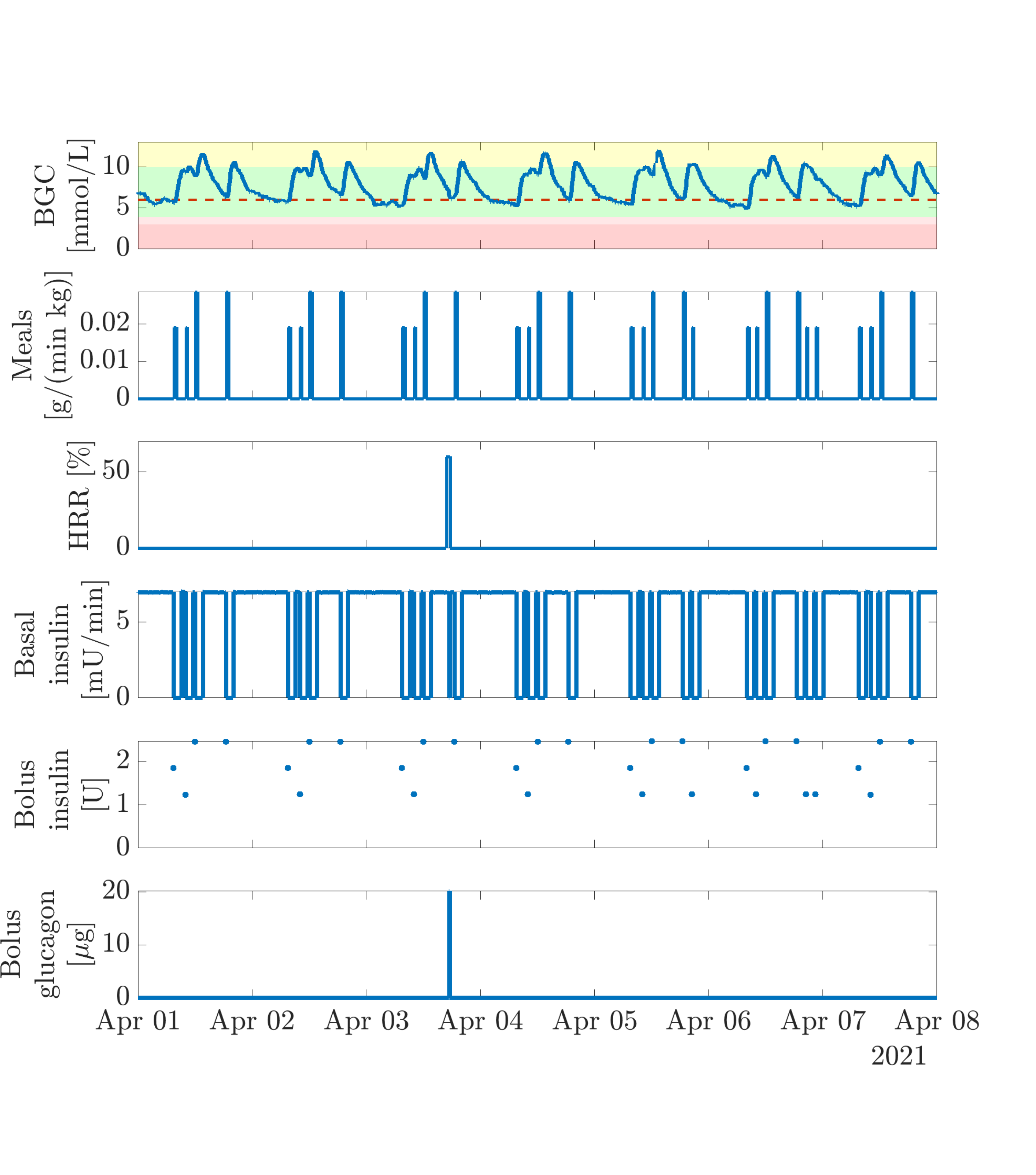}
	\caption{Week 14 of a simulation for 1 patient. From the top we show: 1) The BG concentration (BGC), 2) the meal carbohydrate content (shown as a rate dependent on the body weight), 3) the resting heart rate reserve (HRR), 4) the administered insulin basal rate, 5) the insulin boli, and 6) the glucagon boli. }
	\label{fig:MonteCarloSim3Days210910f5}
\end{figure}

Fig.~\ref{fig:MonteCarloSim3Days210910f5} shows week 14 of a simulation for a single patient. The AP administers insulin boli at meal times, the basal rate is turned off for a period after the meals, and a small glucagon bolus is administered during exercise, i.e., it is not the 100~$\mu$g bolus at the beginning of exercise. As is evident, the TIR is high for the shown period.

\subsection{Performance of the artificial pancreas}\label{sec:performance}
\begin{figure}[tb]
\centering
\includegraphics[trim=5 5 50 5, clip,width=0.45\textwidth]{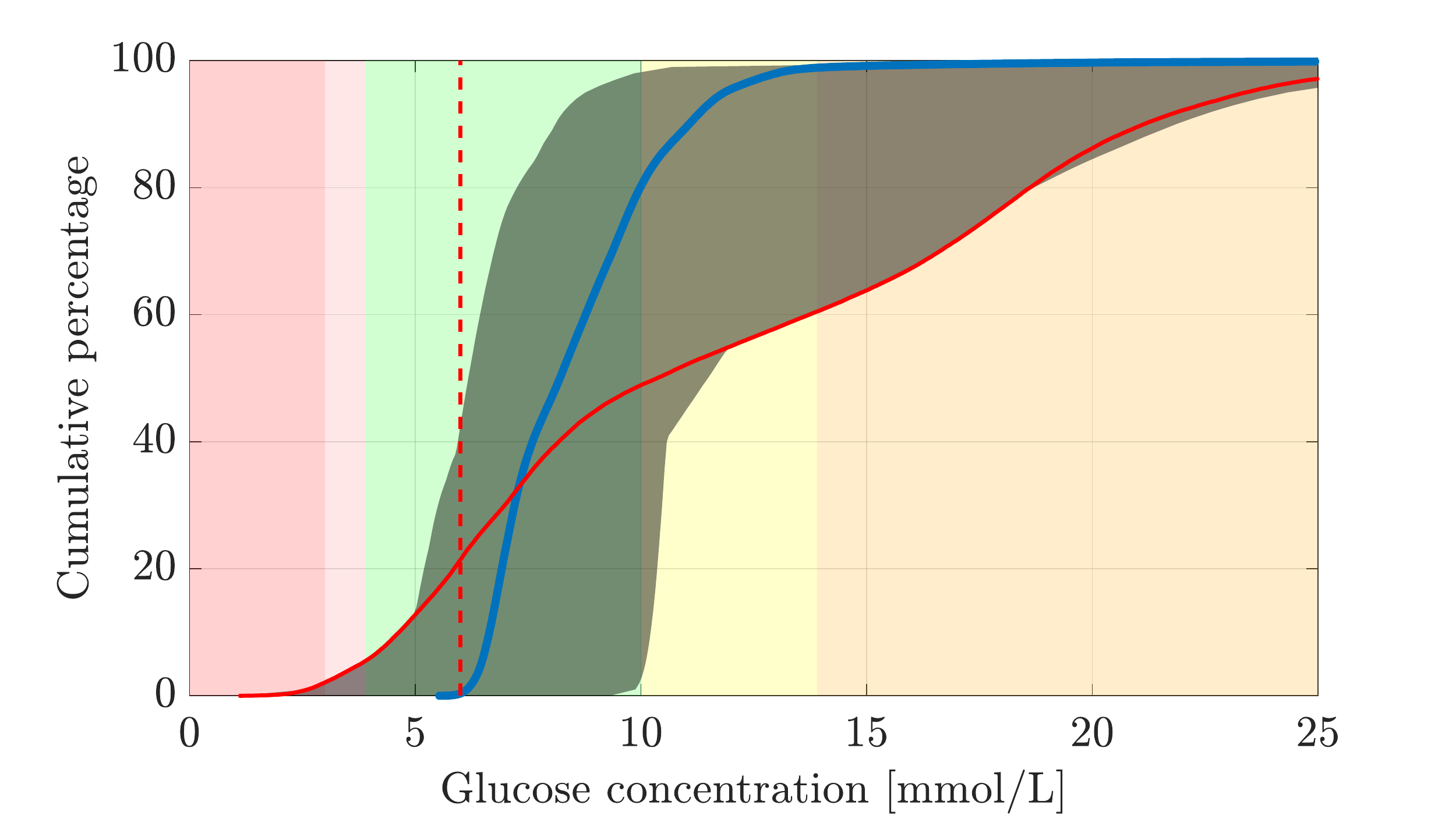}
\caption{Cumulative distribution of the BG concentration in a virtual clinical trial with 1 million patients over 52 weeks. Blue solid line: The mean BG concentration. Red solid line: The worst-case patient. Red dashed line: The setpoint. Grey shaded area: The span of all the patients.}
\label{fig:MonteCarloSimLowMemoryf3}
\end{figure}

\begin{figure*}[tb]
\centering
\includegraphics[trim=0 40 0 30, clip,width=0.175\textwidth]{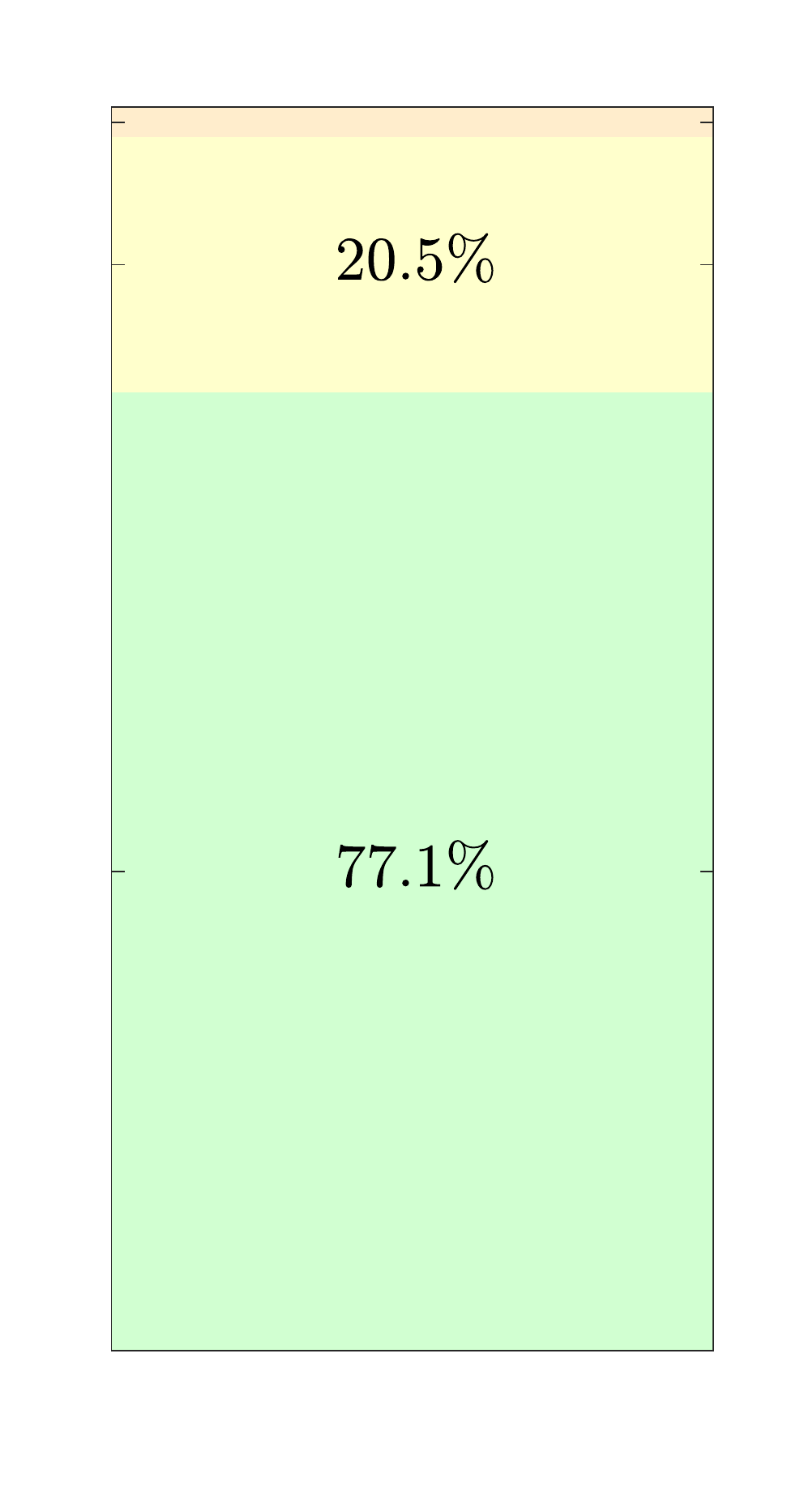}~
\includegraphics[trim=0 40 0 30, clip,width=0.175\textwidth]{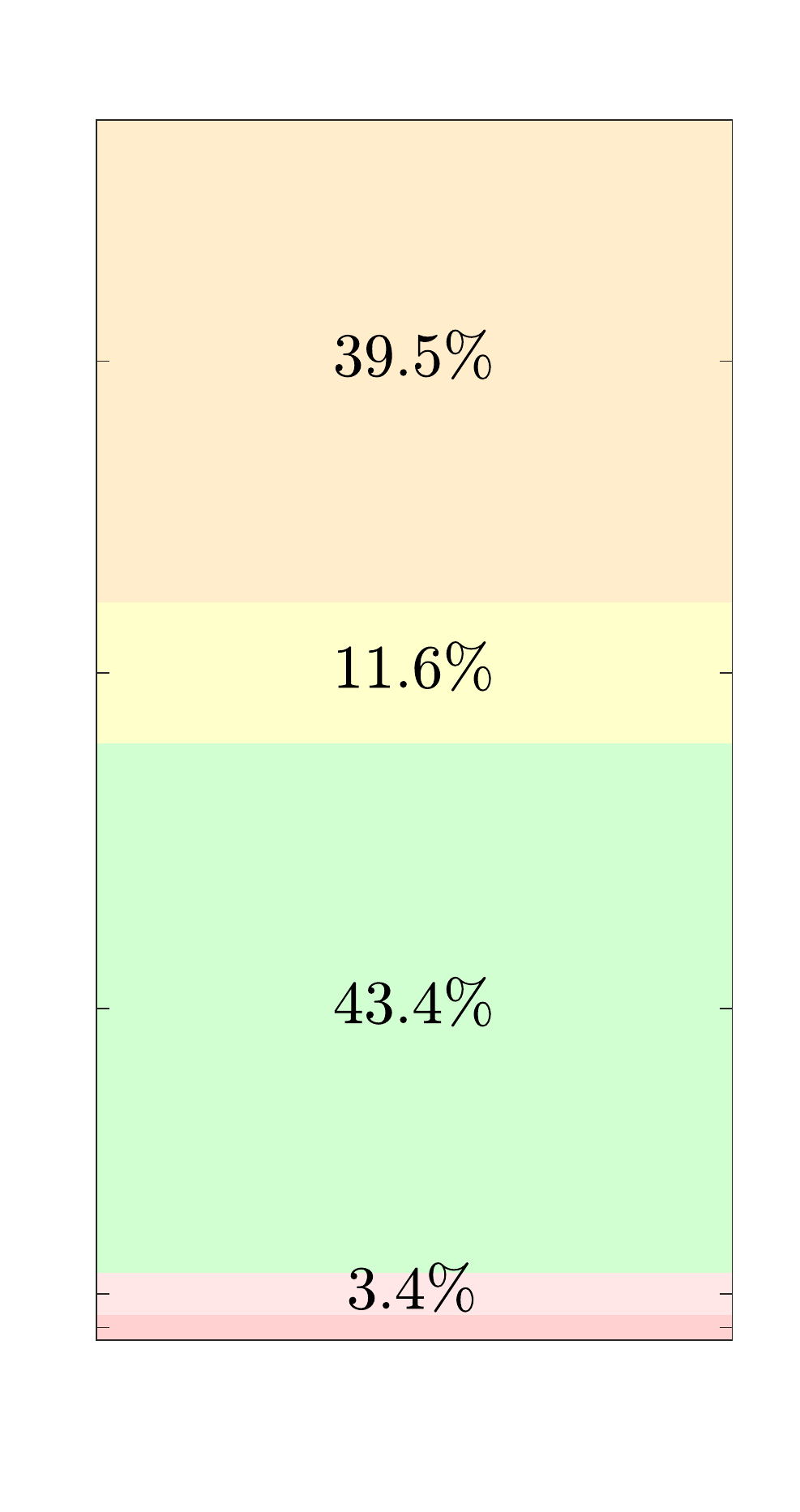}~
\includegraphics[trim=0 34 10 20, clip,width=0.57\textwidth]{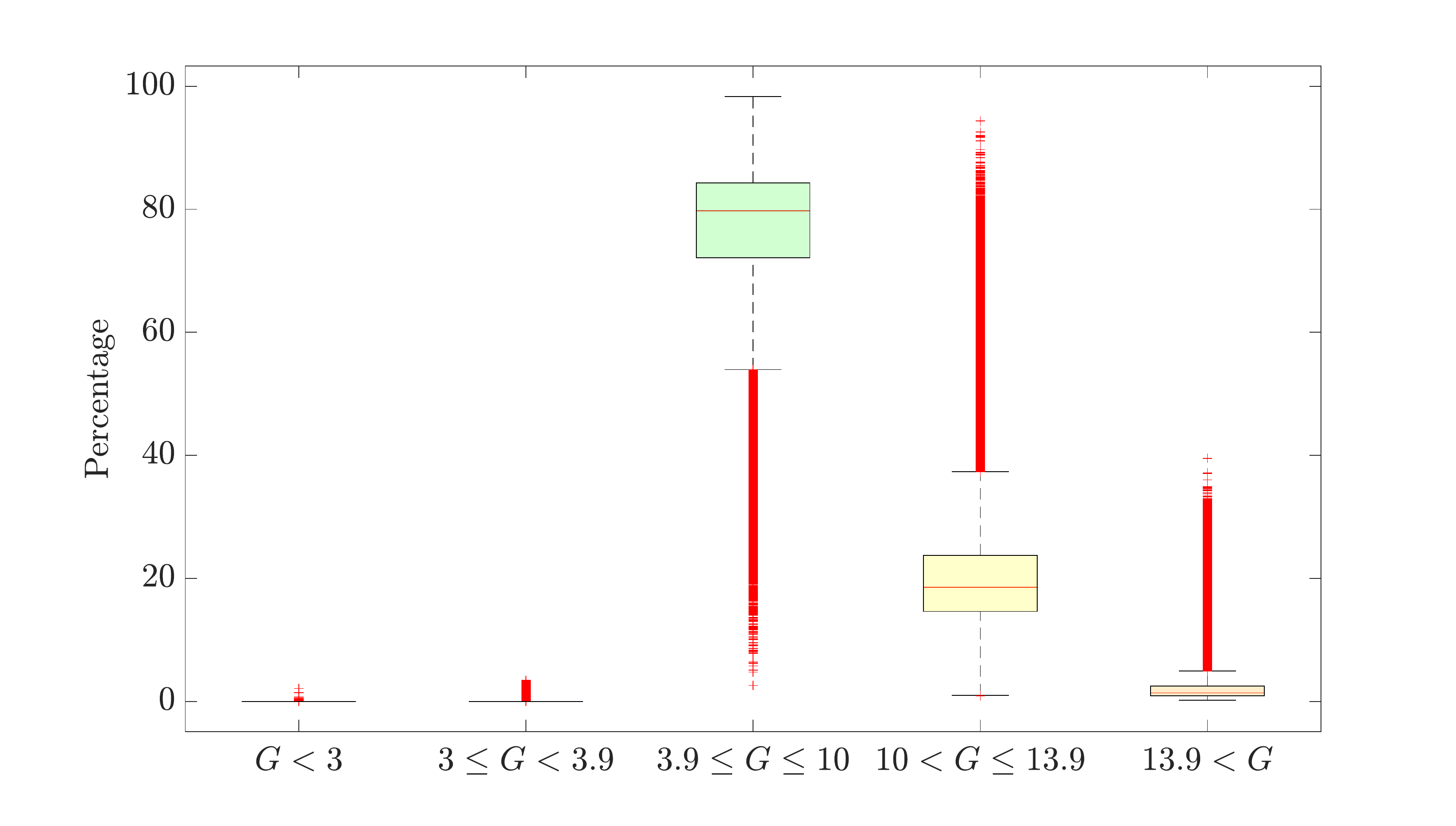}
\caption{Distribution of the TIRs for a virtual clinical trial with 1 million patients over 52 weeks. Left: Mean TIRs for all patients. Middle: TIRs for the worst-case patient. Right: Box plot of the TIRs for all patients.}
\label{fig:MonteCarloSimLowMemoryf1}
\end{figure*}

Fig.~\ref{fig:MonteCarloSimLowMemoryf3} shows the amount of time spent below different BG concentrations. It allows us to inspect the worst-case patient (lowest BG concentration reached), the population average, and the values reached by at least one patient (the span). The worst-case patient is useful for identifying weaknesses in the AP. Here, it seems unlikely that the basal rate is too high because the patient suffers from severe hyperglycemia almost 40\% of the time. This can also be seen from the stacked bar chart in the middle of Fig.~\ref{fig:MonteCarloSimLowMemoryf1}. The population mean can reveal systemic issues. On average, the patients spend almost all of their time above the target of 6~mmol/L. Perhaps too little insulin is administered. Finally, the shaded area can be used to conclude that, e.g., 1)~no patient spends more than 6\% of their time below 3.9~mmol/L and 2)~nobody is above 13.9~mmol/L more than 40\% of the time.

\begin{figure}[tb]
\centering
\includegraphics[trim=10 45 10 45, clip,width=0.45\textwidth]{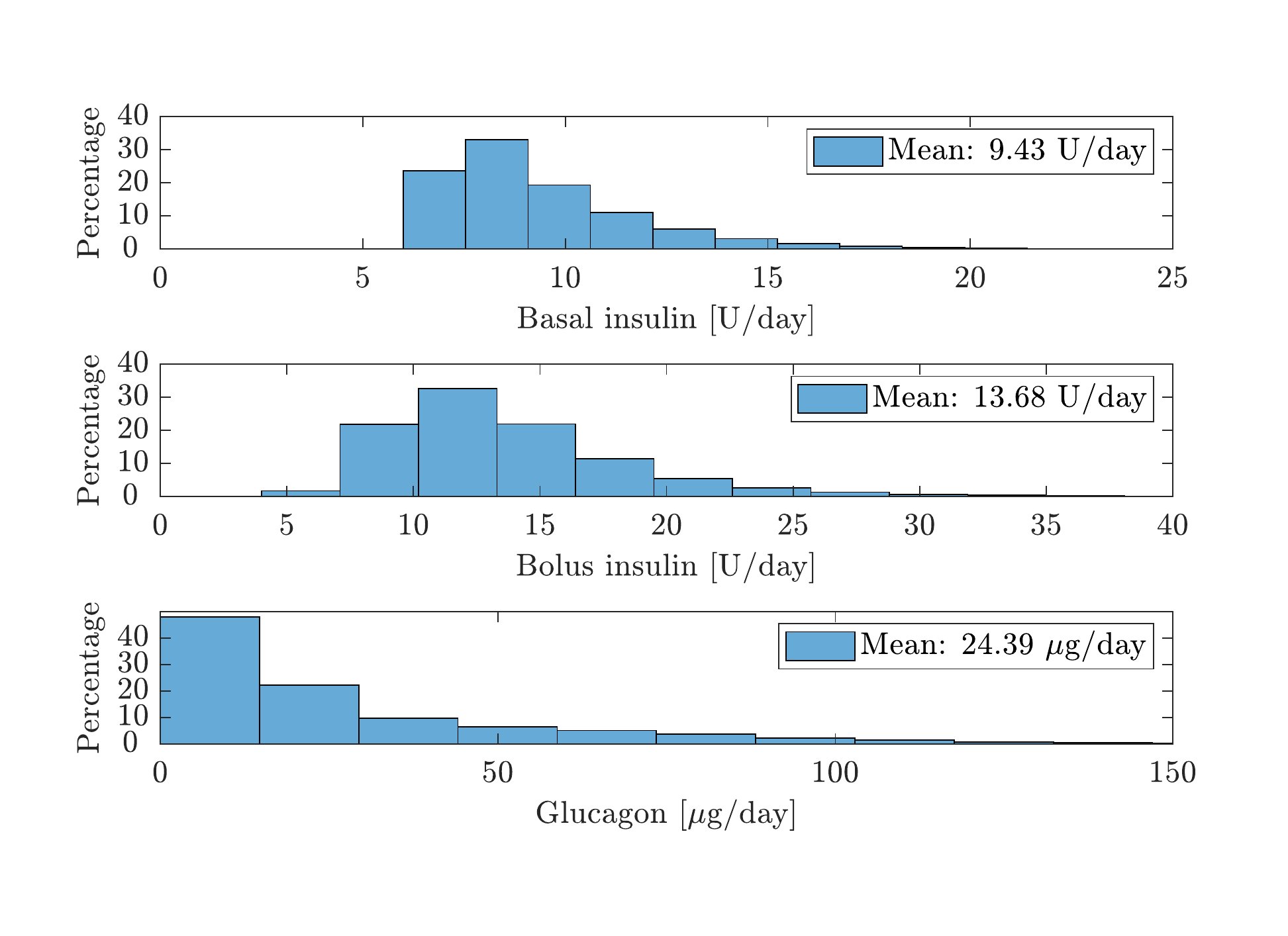}
\caption{Distributions of the total daily basal and bolus insulin and bolus glucagon for a virtual clinical trial with 1 million patients over 52 weeks. Top: Basal insulin. Middle: Bolus insulin. Bottom: Bolus glucagon.}
\label{fig:MonteCarloSimLowMemoryf6}
\end{figure}

Whereas the two stacked bar charts in the left of Fig.~\ref{fig:MonteCarloSimLowMemoryf1} provide an intuitive overview, the box plot on the right gives a comprehensive picture of the TIR for the entire population. The red markers show that only a few patients experience hypoglycemia. On average, the patients spend 77\% of their time in range, and for most patients, this value is at least 55\%. However, a significant part spends between 5\% and 40\% above 13.9~mmol/L. Finally, Fig.~\ref{fig:MonteCarloSimLowMemoryf6} shows distributions of the total insulin and glucagon administered per day. These distributions can be compared with dosage guidelines to see if extreme amounts are administered. For instance, it is positive that for the most part, only small amounts of glucagon are administered here.

\subsection{Comparison of scenarios}\label{sec:comparison}
\begin{figure}[tb]
	\centering
	\includegraphics[trim=5 5 50 5, clip,width=0.45\textwidth]{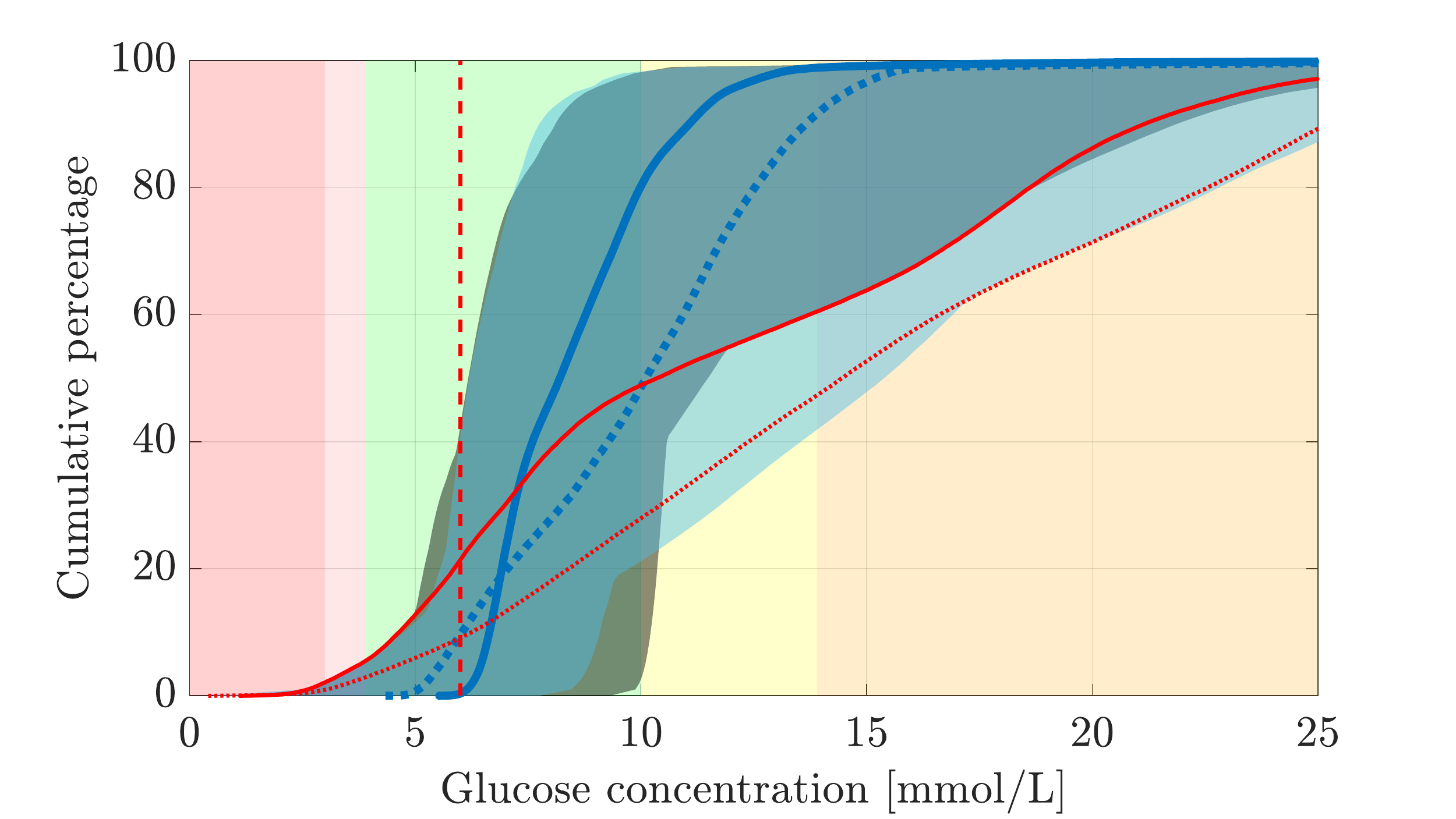}
	\caption{Cumulative distribution of the BG concentration in two virtual clinical trials of 1 million patients over 52 weeks. Blue solid line: The mean BG concentration of trial A. Blue dotted line: The mean BG concentration of trial B. Red solid line: The patient that reaches the lowest BG concentration in trial A. Red dotted line: The patient that reaches the lowest BG concentration in trial B. Red dashed line: The setpoint. Grey shaded area: The span of all patients in trial A. Light blue shaded area: The span of all patients in trial B.}
	\label{fig:MonteCarloSim2Controllers7Daysf3}
\end{figure}

\begin{figure*}[tb]
	\centering
	\includegraphics[trim=0 30 0 30, clip,width=0.175\textwidth]{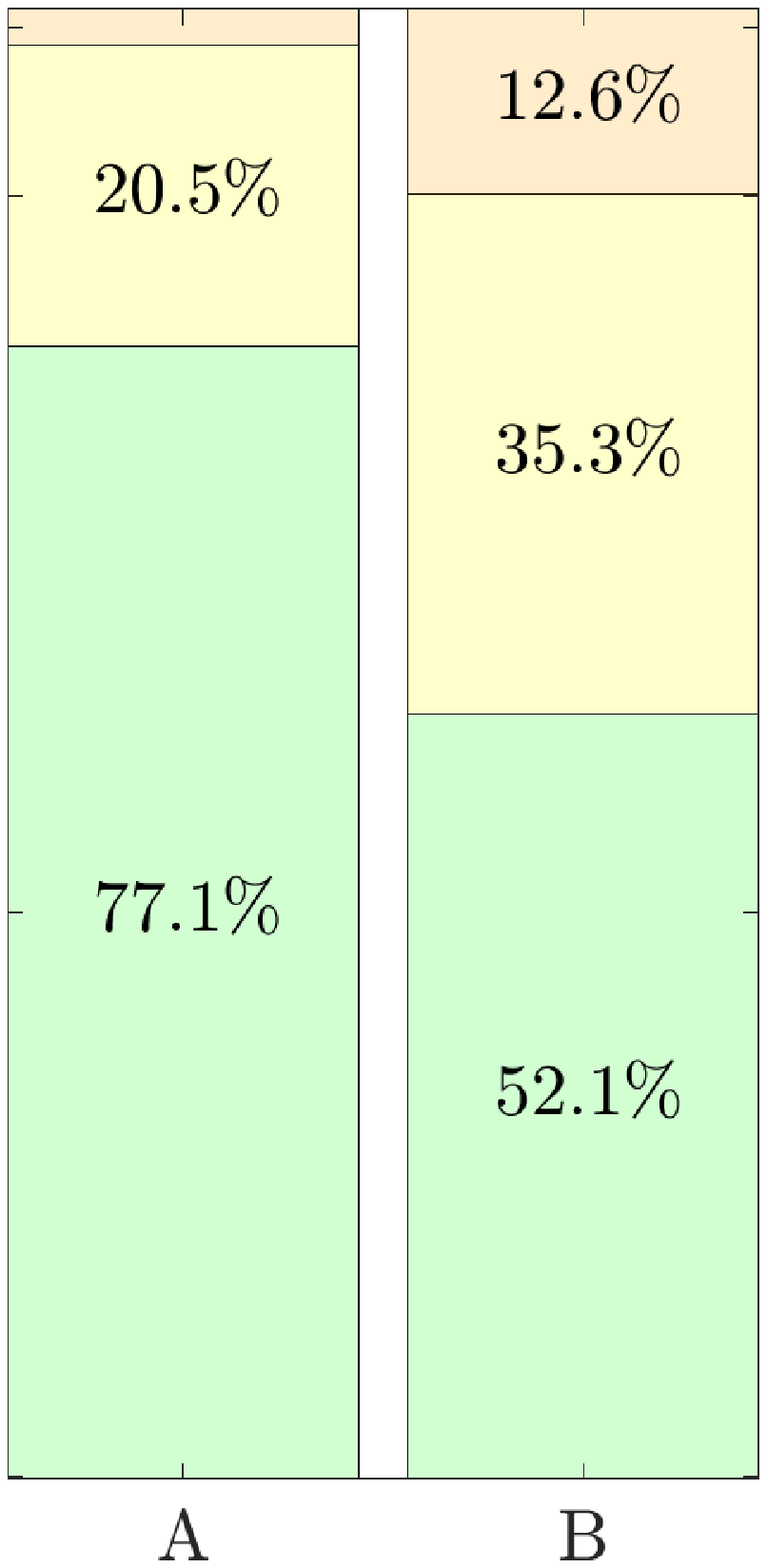}~
	\includegraphics[trim=0 30 0 30, clip,width=0.175\textwidth]{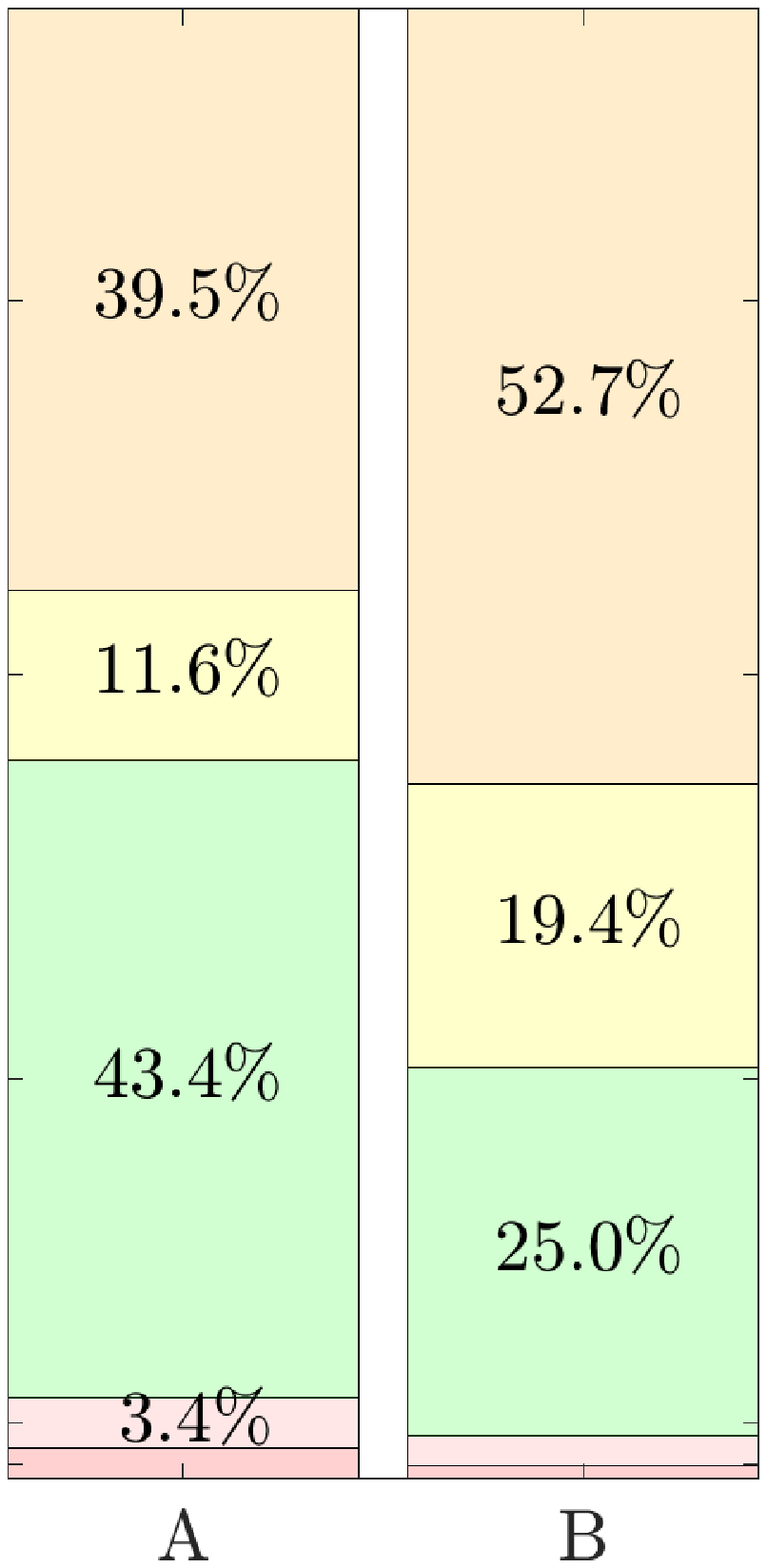}~
	\includegraphics[trim=0 30 10 10, clip,width=0.57\textwidth]{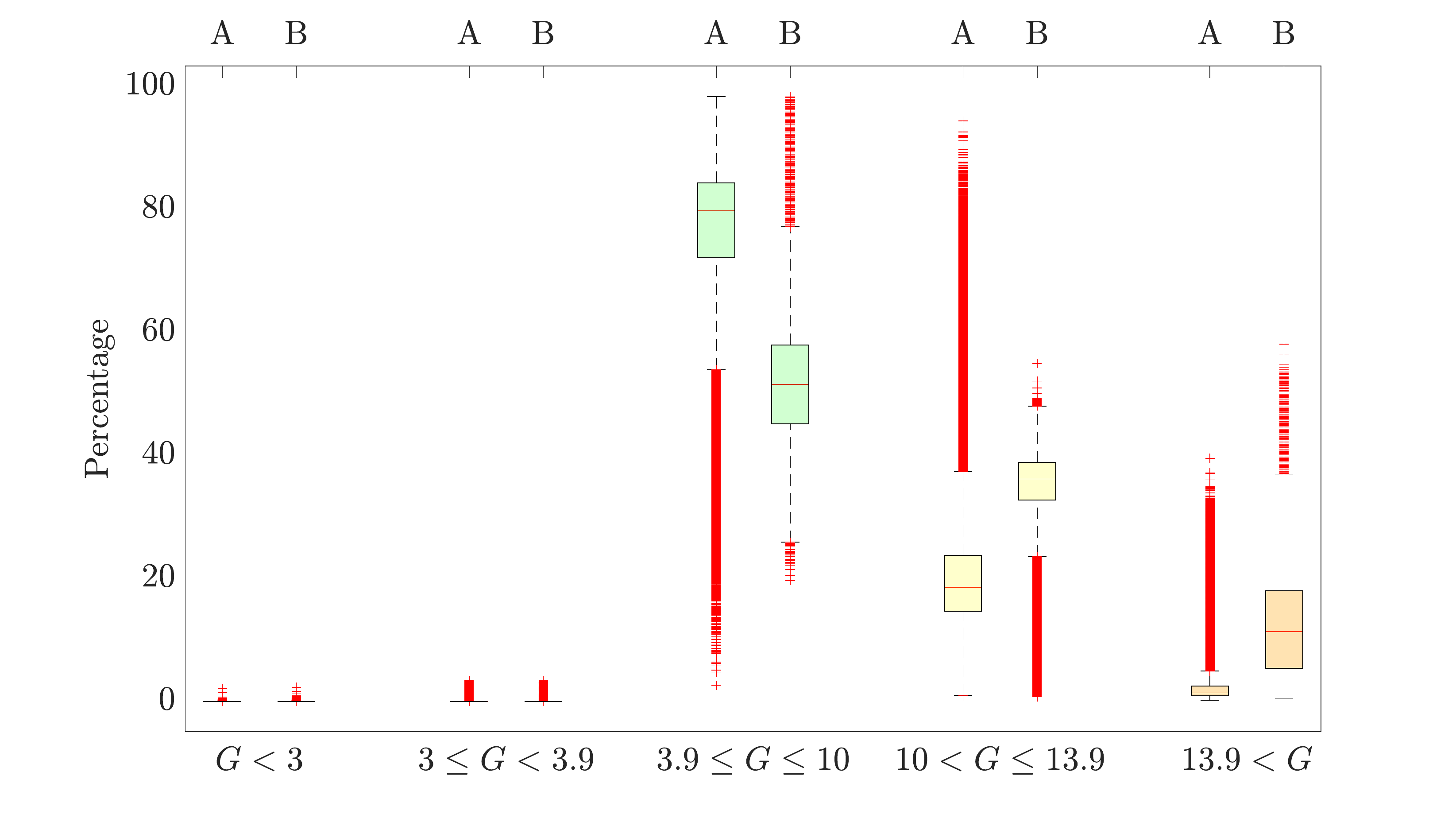}
	\caption{Distribution of the TIRs for two virtual clinical trials of 1 million patients over 52 weeks. Left: Mean TIRs for all patients in trial A and trial B. Middle: TIRs for the worst-case patients in trial A and trial B. Right: Box plot of the TIRs for all patients in trial A and trial B.}
	\label{fig:MonteCarloSim2Controllers7Daysf4}
\end{figure*}

\begin{figure}[tb]
	\centering
	\includegraphics[trim=10 45 10 45, clip,width=0.45\textwidth]{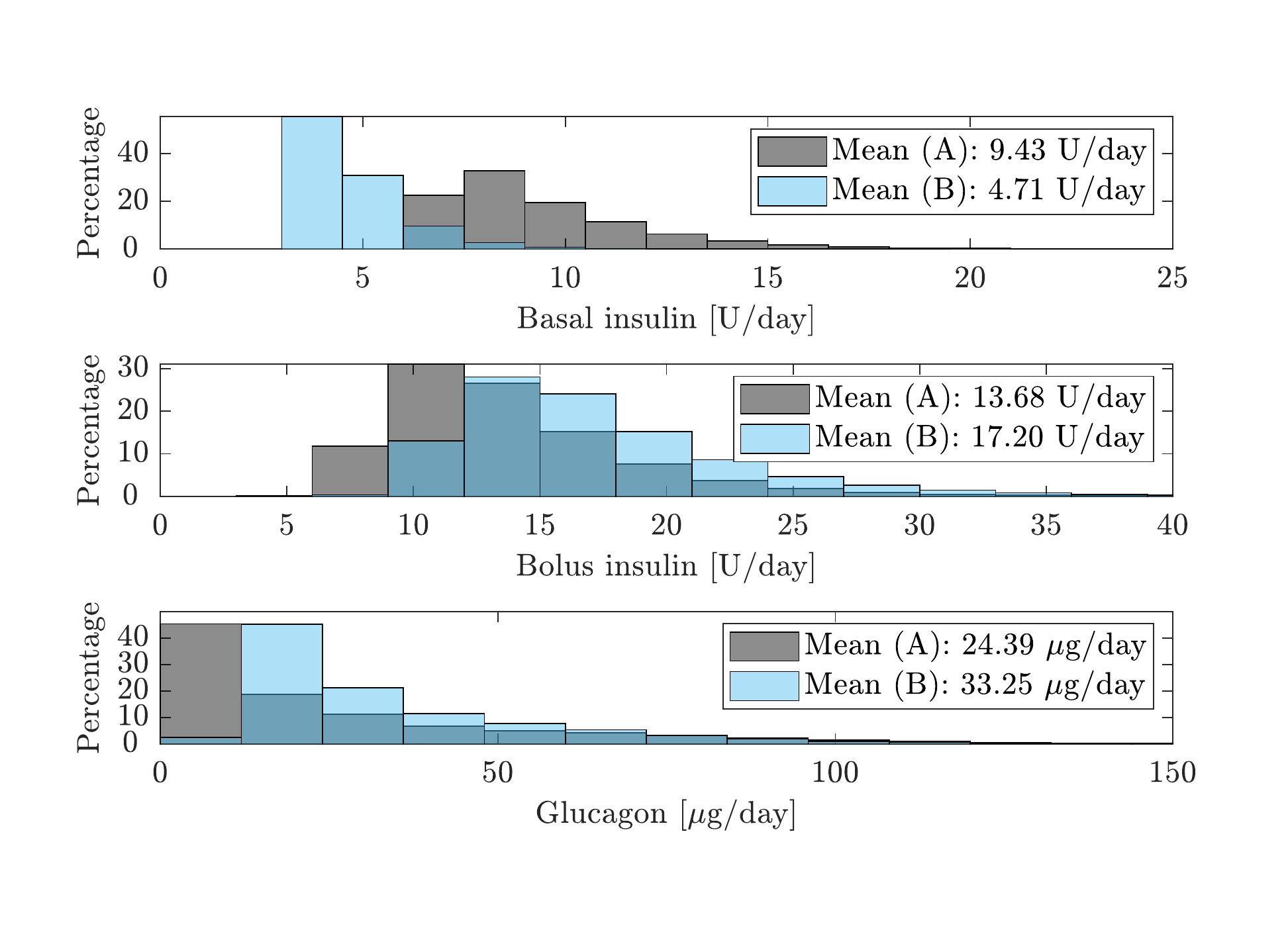}
	\caption{Distributions of the total daily basal and bolus insulin and bolus glucagon for two virtual clinical trials with 1 million patients over 52 weeks. Top: Basal insulin. Middle: Bolus insulin. Bottom: Bolus glucagon. Grey: Trial A. Light blue: Trial B.}
	\label{fig:MonteCarloSim2Controllers7Daysf6}
\end{figure}

When comparing different scenarios or APs, the plots in Fig.~\ref{fig:MonteCarloSimLowMemoryf3}--\ref{fig:MonteCarloSimLowMemoryf6} can be overlaid or combined as shown in Fig.~\ref{fig:MonteCarloSim2Controllers7Daysf3}--\ref{fig:MonteCarloSim2Controllers7Daysf6}. It is clear from Fig.~\ref{fig:MonteCarloSim2Controllers7Daysf3} that the basal rate is too low in trial~B. However, the worst-case patients reach equally low BG concentrations in both trials. Fig.~\ref{fig:MonteCarloSim2Controllers7Daysf4} also clearly shows that far better TIRs are reached with the correct basal rate although the worst-case patient experiences hypoglycemia more often with the correct basal rate. Fig.~\ref{fig:MonteCarloSim2Controllers7Daysf6} allows for direct comparison of how much insulin and glucagon that is administered on a daily basis. Obviously, less basal insulin is administered in trial B, and as a result, more bolus insulin is given. Consequently, more glucagon is given in trial B.


\section{Conclusion}\label{sec:Conclusion}
In this paper, we present a virtual clinical trial for assessing the uncertainty in the performance of closed-loop diabetes treatments. We use a high-performance closed-loop Monte Carlo method for quantifying the uncertainty, and we evaluate the performance by examining 1)~the distributions of the TIRs for the patient population and 2)~the distributions of the total daily doses of basal and bolus insulin as well as bolus glucagon.
Furthermore, this approach can be used to compare the performance in different scenarios and for different closed-loop treatments.
Finally, we demonstrate that a virtual clinical trial with one million patients over 52 weeks can be completed in 1~h and 22~min by using parallel high-performance software and hardware. The developed software can be used for closed-loop systems with any drug administration device (pump or pen) and measurement device (CGM or self-monitoring of blood glucose (SMBG) device).

\bibliographystyle{IEEEtran}
\bibliography{bib/IEEEabrv,ref/References}
\end{document}